%% file: main_1.tex
\documentclass[review]{elsarticle}

\usepackage{lineno, hyperref, siunitx}
\usepackage{amssymb,amsmath}
\usepackage[margin=1in]{geometry}
\usepackage{textpos}

\usepackage{enumitem}
\usepackage[parfill]{parskip}    
\usepackage{latexsym, amssymb, bbm}
\usepackage{amsmath}
\usepackage{amsfonts}
\usepackage{amsthm}
\usepackage{verbatim}
\usepackage{fancyhdr}
\usepackage{mathrsfs}
\usepackage{enumerate}
\usepackage{geometry}
\usepackage{graphicx}
\usepackage{algorithm}
\usepackage{algpseudocode}
\usepackage{caption}
\usepackage{subcaption}
\usepackage{cancel}
\usepackage{listings}
\usepackage{xcolor}
\usepackage{booktabs}
\usepackage{listings}
\modulolinenumbers[5]

\journal{arXiv}

\begin{document}

\begin{frontmatter}

\title{Ensemble WSINDy for Data Driven Discovery of Governing Equations \\ from Laser-based Full-field Measurements}

\author[cvenaddress]{Abigail C. Schmid}
\author[asenaddress, appmaddress]{Alireza Doostan\corref{alireza}}
\cortext[alireza]{alireza.doostan@colorado.edu}
\author[cvenaddress, appmaddress]{Fatemeh Pourahmadian\corref{fatemeh}}
\cortext[fatemeh]{fatemeh.pourahmadian@colorado.edu}

\address[cvenaddress]{Department of Civil, Environmental and Architectural Engineering, University of Colorado Boulder}
\address[asenaddress]{Ann \& H.J. Smead Department of Aerospace Engineering Sciences, University of Colorado Boulder}
\address[appmaddress]{Department of Applied Mathematics, University of Colorado Boulder, USA}

\begin{abstract}

    This work leverages laser vibrometry and the weak form of the sparse identification of nonlinear dynamics (WSINDy) for partial differential equations to learn macroscale governing equations from full-field experimental data. In the experiments, two beam-like specimens, one aluminum and one IDOX/Estane composite, are subjected to shear wave excitation in the low frequency regime and the response is measured in the form of particle velocity on the specimen surface. The WSINDy for PDEs algorithm is applied to the resulting spatio-temporal data to discover the effective dynamics of the specimens from a family of potential PDEs. The discovered PDE is of the recognizable Euler-Bernoulli beam model form, from which the Young's modulus for the two materials are estimated. An ensemble version of the WSINDy algorithm is also used which results in information about the uncertainty in the PDE coefficients and Young's moduli. The discovered PDEs are also simulated with a finite element code to compare against the experimental data with reasonable accuracy. Using full-field experimental data and WSINDy together is a powerful non-destructive approach for learning unknown governing equations and gaining insights about mechanical systems in the dynamic regime.  
    
\end{abstract}

\begin{keyword}
    equation discovery; WSINDy; data driven mechanics; ensemble methods; full-field measurements; laser vibrometry 
\end{keyword}

\end{frontmatter}

\section{Introduction}
\label{sec:intro}

    In science and engineering, differential equations are powerful tools for describing, analyzing, and making predictions about a system of interest. These equations have historically been derived from theory and first principles. However, many systems are too complex or lacking information to be fully derived in this way, thus making the governing differential equations unknown or only partially known. In addition, accurate resolution of a complex PDE system may become computationally prohibitive. These pose challenges for having useful models of the system. 
    
    One area where this difficulty is apparent is in studying the mechanics of composite particulate materials. Understanding the behavior of particulate materials is generally a multiscale problem due to the interaction of the particle (micro) and bulk material (macro) scales. Even only considering the macroscale, modeling these materials can be challenging since the governing equations may be unknown and bulk level material properties can be difficult to obtain from the various materials and likely random microstructure. An interesting class of composite particulate materials is plastic bonded explosives (PBXs), which are comprised of a high explosive (HE) crystal and a polymeric binder. It has been found that different features of PBXs, such as density, grain size distribution of the crystals, manufacturing process, and temperature, change the performance and mechanical properties of the material \cite{HermanComposite2021, LiuPBX95012020}. It is important to understand these effects on the mechanical behavior of the HE for modeling and safety purposes. However, performing mechanical testing on HEs is dangerous and costly, thus many aspects of these observations have not yet been characterized. One approach to better understand HEs involves the use of mock HE materials, or inert materials which mimic properties and some behavior of the HE without the risk of detonation. For example, a combination of IDOX (5-iodo-2'-deoxyuridine) crystals and nitroplasticized Estane 5703 binder has recently been formulated as a mock for the HE PBX 9501, which uses the same binder \cite{HermanComposite2021, LiuPBX95012020, YeagerDevelopment2019}. Experimental tests can more safely be conducted on the IDOX/Estane mock and any insights gained from tests on the mock should extend to the HE as well. Mechanical testing has been conducted on the bulk mock composite and the individual components via quasi-static compression \cite{HermanComposite2021, CadyMechanical2006}, Split Hopkinson Pressure Bar \cite{CadyMechanical2006, BurchCompressive2022}, Brazilian disk compression \cite{LiuPBX95012020, YeagerDevelopment2019}, and nanoindentation \cite{YeagerDevelopment2019, BurchNanoindentation2017}. The compression tests result in stress-strain curves from which quantities such as the maximum stress, strain at maximum stress, and macroscale Young's modulus can be determined. These quantities can then be used in modeling the behavior of the mock. However, a more complete description of the bulk behavior such as a full governing equation in the form of a partial differential equation (PDE) would be more informative and useful to extend any insights to HEs. 

    In order to generate a PDE description of the behavior of the IDOX/Estane composite, this study leverages data from laser-based vibrometry tests. In these tests, a 1D sample of the material of interest is excited by shear waves and the resulting wave motion in the sample is measured in terms of transverse particle velocity on the sample surface. Since the macroscale response is the main interest, low frequency input waves are used so that the wavelength is much longer than the cross section of the samples and any microstructure in the composite. Compared to the testing methods listed above which have been used with the IDOX/Estane material \cite{HermanComposite2021, LiuPBX95012020, YeagerDevelopment2019, CadyMechanical2006, BurchCompressive2022, BurchNanoindentation2017}, this method is non-destructive and can be repeated on the same specimen with different inputs. In addition to testing the IDOX/Estane composite, aluminum is used as a test case material since it has been studied more thoroughly. These tests result in a dense spatio-temporal grid of velocity measurements which can be used to generate insights about the mechanics of the material, with appropriate data analysis methods. 

    For a known or assumed governing equation form, one could use the full-field experimental data in the context of an inverse problem to estimate the mechanical properties of the sample. For example, using modal analysis \cite{HeModal2004}, Kalman filters \cite{DelalleauUse2008, NguyenApplication2023}, or physics informed neural networks (PINNs) \cite{RaissiPhysics2019, XuDeep2023, NormanConstrained2024}. Experimental data are inherently noisy, which can pose challenges for generating accurate inversion results. Moreover, since the IDOX/Estane mock is not well characterized, it is unclear what the governing equation should be, so imposing a model form is likely to introduce error. Therefore, this work approaches the problem as one of model identification or equation discovery. 
    
    Data driven methods for equation discovery, such as symbolic regression \cite{SchmidtDistilling2009}, NARMAX \cite{ChenRepresentation1989, WeiNARMAX2005, RahroohIdentification2009}, neural networks \cite{NormanConstrained2024}, and sparse regression methods \cite{BruntonDiscovering2016, SchaefferLearning2017, MessengerGalerkin2021, CortiellaPriori2023, WentzDerivative2023} have the potential to identify unknown or partially known governing equations from full-field data. Focusing on sparse regression type methods, the sparse identification of nonlinear dynamics (SINDy) has become a popular framework for equation discovery for ordinary and partial differential equations \cite{BruntonDiscovering2016, WentzDerivative2023, RudyData2017}. Briefly, in the SINDy framework, a library of potential terms in the differential equation is defined and the data derivatives are computed based on those terms. Then, a sparsity promoting regression scheme is used to identify which terms are active in the data and the corresponding coefficients. This approach accurately recovers models from simulated data but struggles in the presence of measurement noise \cite{CortiellaPriori2023, WentzDerivative2023, HokansonSimultaneous2023}. Due to the challenges of discovering models from noisy data, the use of this approach with experimental data has been limited and primarily for biological systems \cite{HoffmannReactive2019, KahemanLearning2019, SandozSINDy2023, BrummerData2023, ProkopBiological2024}. To improve the robustness to noise, using the weak form of PDEs has also been introduced \cite{MessengerGalerkin2021, MessengerLearning2022, ReinboldUsing2020, GurevichRobust2019}. While this change does improve the robustness to noise, it has primarily been demonstrated via simulated data with added noise \cite{MessengerLearning2022, MessengerOnline2022} and its use with experimental data is, so far, still limited \cite{WoodallModel2024}.

    In this work, the weak form of the sparse identification of nonlinear dynamics (WSINDy) for PDEs is the chosen approach for the equation discovery task with the experimental data. Briefly, the WSINDy approach \cite{ReinboldUsing2020, GurevichRobust2019, MessengerWeak2021, WangVariational2019} performs sparse regression on a library of potential PDE terms to identify which of the terms are active in the data and the coefficients on those terms. The WSINDy approach has several desirable features which make it an appropriate choice for this application: (1) the weak form makes the method more robust to noise in experimental measurements, (2) the method discovers PDEs resulting in a governing equation for the spatio-temporal system, (3) the method is boundary condition agnostic due to the construction of the test functions, thus avoiding the need to estimate boundary conditions from the experimental data. This study utilizes the WSINDy algorithm implementation given in \cite{MessengerWeak2021} and presents an extension of the implementation which performs an ensemble version of the algorithm using subsamples of the data. The goal is to recover uncertainty information on the terms and coefficients in the discovered equations. Other ensemble SINDy approaches have focused on synthetic data, bagging the library terms, and bootstrapping the data \cite{FaselEnsemble2022}.  
    
    Overall, this work aims to learn macroscale governing equations, with quantified uncertainty, from the experimental data. The main contributions of this study are three-fold. First, the WSINDy framework is successfully applied to discover governing equations from full-field experimental data. Here it is shown that the method works to discover interpretable equations from experimental measurements and provides an example of using the model coefficients to recover a mechanical property of interest from the resulting equations. Additionally, the discovered equations are used with finite element simulations to compare against the experimental data. Second, analyses of these data are inspired by elastography, a technique for reconstructing the elastic modulus from full-field data that was originally designed for medical imaging of soft tissues \cite{MuthupillaiMagnetic1995, BarboneElastic2004}. Typically, the analyses are conducted in the frequency domain \cite{XuDeep2023, PourahmadianElastic2018}; instead, this study is an example of analyzing such data in the time domain. Last, there are only a few data sets pertaining to the IDOX/Estane mock, so the laser vibrometry test data and resulting governing equations provide new information and insights about the material and its macroscopic property. 
    
    The rest of the paper is organized as follows. Sections~\ref{sec:expmethods} and ~\ref{sec:preprocess} describe the experimental procedure and pre-processing of experimental data. Section~\ref{sec:eqdiscoverymethods} outlines the equation discovery and ensemble discovery methods. Section~\ref{sec:results} provides results of the equation discovery and comparison with finite element simulations of the discovered equations. Lastly, Section~\ref{sec:conclusion} summarizes the findings of this work. 
    
\section{Experimental Campaign}
\label{sec:expmethods}

\begin{figure}[h]
        \centering
        \includegraphics[width=0.88\textwidth]{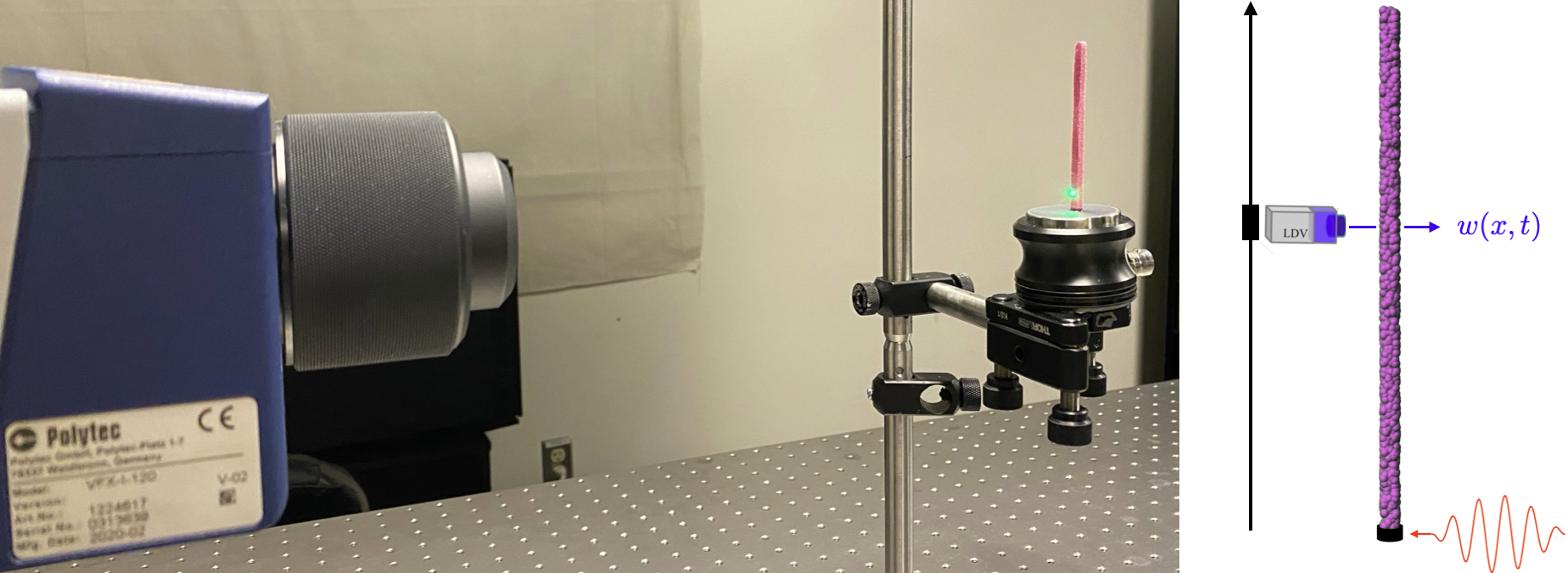}
        \caption{Experimental set-up and sensing configuration diagram. The IE sample is vertically mounted on a shear-wave transducer inducing the five-cycle burst input. The resulting wave motion is then measured in the form of transverse particle velocity by a laser Doppler vibrometer scanning on a uniform grid of $113$ points, with $0.5$ \unit{\milli\meter} increments, along the beam length.}
        \label{fig:expsetup}
    \end{figure}

    Experiments are conducted on two beam-like specimens of aluminum (Al) and a composite of IDOX and Estane (IE). The Al specimen is a cylinder with diameter $6.35$ \unit{\milli\meter} and measurement domain length $97$ \unit{\milli\meter}. The IE specimen is a rectangular prism with cross section $4.18$ $\!\times\!$ $2.84$ \unit{\milli\meter\cubed} and measurement domain length $56.0$ \unit{\milli\meter}, manufactured according to the procedure in~\cite{CamerloPressing2021}. The mass densities are measured at $2721.9$ and $1301.4$ \unit{\kilo\gram\per\meter\cubed} for Al and IE, respectively. The testing configuration, in terms of excitation and sensing, is the same for both specimens as shown in Fig.~\ref{fig:expsetup}. A piezoelectric shear-wave transducer (Olympus V150-RB) attached to the sample's base generates a 5-cycle burst of the form 
    \begin{equation} 
     H(f_c t) \sin(0.2\pi f_c t) \sin(2\pi f_c t) H(5-f_c t), 
    \end{equation}
    
    where $H$ is the Heaviside function and the center frequency $f_c$ is $10$ \unit{\kilo\hertz} for the Al test and $3$ \unit{\kilo\hertz} for the IE test. These input frequencies represent the long wavelength regime for the respective specimens. The induced wave motion is then measured in terms of particle velocity on the sample surface using the Polytec VibroFlex (Xtra VFX-I-120) Laser Doppler Vibrometer (LDV) which is mounted on a programmable robotic arm for scanning to perform equally spaced measurements along the length of the specimens. For both specimens, measurements are taken in $0.5$ \unit{\milli\meter} increments resulting in $195$ measuring points for the Al specimen and $113$ for the IE. At each scan point, the time history of velocity response is captured. For the Al specimen, the time series is measured in increments of $1.600e-08$ \unit{\second} over $0.002$ \unit{\second} and for the IE specimen the data is collected every $3.195e-08$ \unit{\second} over $0.004$ \unit{\second}. The sampling time steps are determined based on the center frequency of the input signal and sampling rates of the equipment. In general, the LDV system is capable of tracking particle motion with magnitudes down to O(\text{nm}) in terms of displacement, spatial resolution of O(100$\mu$m), and frequency range $\sim$DC to 24 MHz~\cite{XuDeep2023, PourahmadianElastic2018}. 

\section{Pre-Processing of Experimental Data}
\label{sec:preprocess}

\begin{figure}[]
        \centering
        \includegraphics[width=0.9\textwidth]{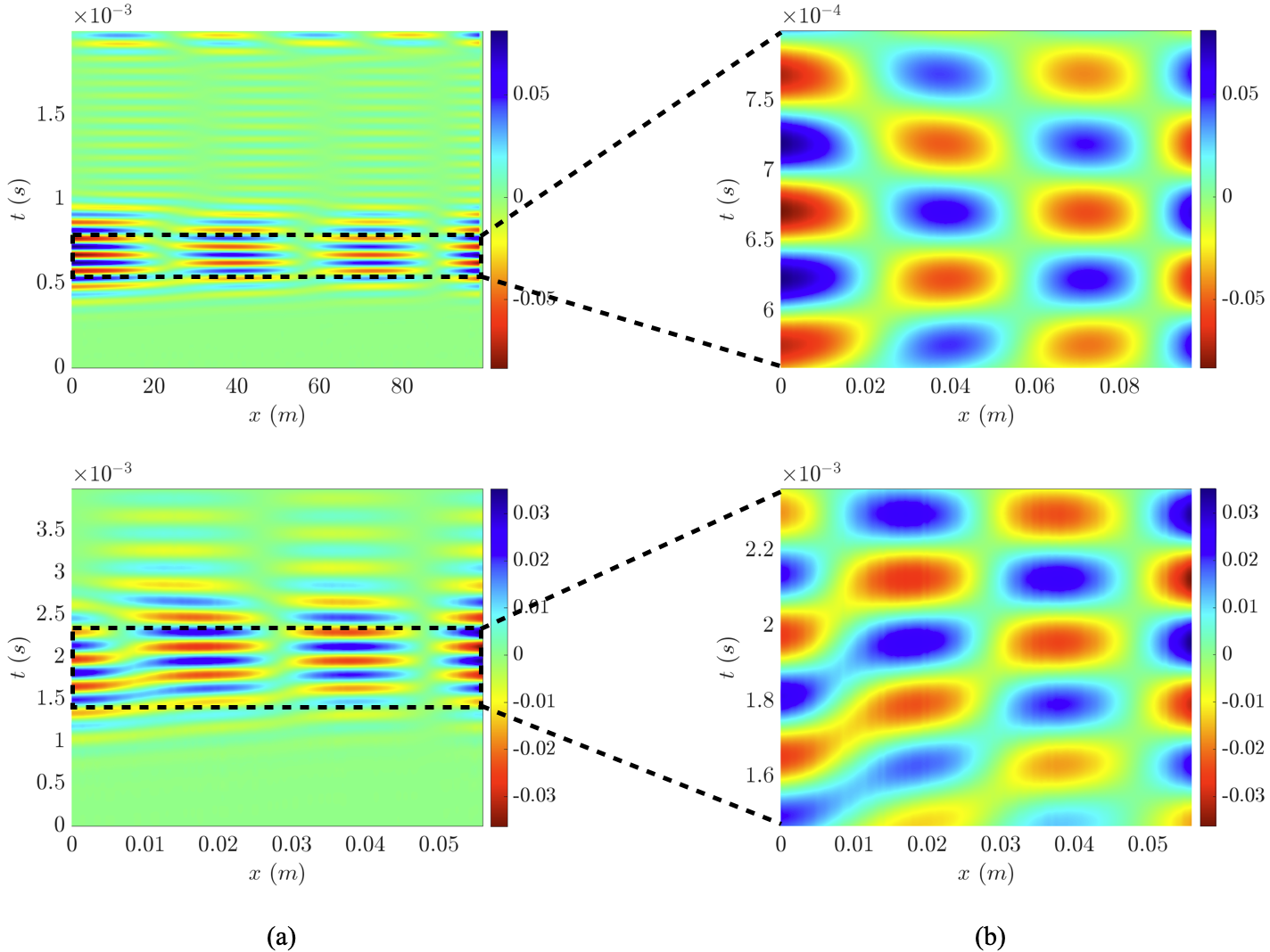}
        \caption{Experimental data. Each panel shows the velocity measurements (colorbar) in time (y-axis) and space (x-axis). The aluminum data is in the top row and the IE is in the bottom row. The filtered experimental data is shown in (a) and the truncated data used for equation discovery is shown in (b).}
        \label{fig:exp-data}
    \end{figure}

    The raw experimental data is pre-processed before being used for equation discovery. First, the data is downsampled in time by a factor of 10 and a bandpass filter is applied to each dataset. The filter removes frequencies outside the range of the input signal which represent noise in the system introduced by the experimental equipment. In the case of the Al test conducted with the input signal centered at $10$ \unit{\kilo\hertz}, the bandpass filter removes frequencies outside of $4$ to $16$ \unit{\kilo\hertz}. For the IE test conducted with the input centered at $3$ \unit{\kilo\hertz}, the bandpass filter removes frequencies outside the range $1$ to $5$ \unit{\kilo\hertz}. The filtered data is shown in Fig.~\ref{fig:exp-data} (a). The entire time series of measurements contains sections with approximately zero motion at the beginning, which is the time between when the input signal is generated and communicated to the transducer, and at the end after the main signal has passed. Since the goal is to learn governing equations for the motion of the specimens, the data is truncated to only include the sections containing the signal with the largest magnitude. For the Al data, this window is between $5.5984e-04$ and $7.9984e-04$ \unit{\second} with $dt = 1.600e-07$ \unit{\second} and for the IE it is between $1.4693e-03$ and $2.3639e-03$ \unit{\second} with $dt = 3.195e-07$ \unit{\second}. These windows are shown in Fig.~\ref{fig:exp-data} (b). 
    
\section{Equation Discovery Methods}
\label{sec:eqdiscoverymethods}

\subsection{Weak SINDy for PDEs}
\label{sec:wsindy}

    \subsubsection{Problem Statement and Weak Form}

    The WSINDy algorithm as implemented in~\cite{MessengerWeak2021} is used to discover partial differential equations (PDEs) from the aforementioned experimental data sets. As will be explained in the following, using notation similar to that used in \cite{MessengerWeak2021}, the WSINDy algorithm uses a weak formulation and sparse regression to discover relevant differential operators and functions governing the behavior of the system from a library of possible equations. 
    
    Let $w(x,t)$ be the field variable as a function of space and time. In this work, $w$ is the flexural velocity. Then, let $\mathbf{W} \in \mathbb{R}^{N_x \times N_t}$ describe a matrix containing measurements of $w$ on a spatio-temporal grid of spatial points $x_n$, $n = 1,\ldots, N_x$, and time samples $t_m$, $m = 1,\ldots, N_t$ such that $W_{nm} = w(x_n, t_m) + \varepsilon_{nm}$ where $\mathbf{\varepsilon} \in \mathbb{R}^{N_x \times N_t}$ is a matrix of independent and identically distributed random noise. Thus, $\mathbf{W}$ collects the experimental velocity data over the space-time grid $(\mathbf{x},\mathbf{t}) \in [0, X]^{N_x}\! \times [0, T]^{N_t}$. 
    
    The WSINDy for PDEs algorithm takes the spatio-temporal data set $\mathbf{W}$ and grids of of time points $\mathbf{t}$ and spatial measurement points $\mathbf{x}$ as inputs. The algorithm also requires that the user defines the term on the left-hand side of the governing equation to be discovered. Here that term was chosen as the second-order time derivative of the field variable, $w_{,tt}$, where the subscripts denote the order and variable of differentiation. Thus, a general strong form of the PDE can be written as
    \begin{equation} w_{,tt} = F(w), \label{eq:wsindy-general-wtt}\end{equation}where $w_{,tt}$ is again the term chosen for the left-hand side, $w$ is the velocity field variable, and $F$ is a nonlinear differential operator to be identified. 

    In order to convert Eq.~\ref{eq:wsindy-general-wtt} to its weak form, WSINDy constructs a real-valued, smooth test function $\psi(x,t)$ which is compactly supported in $[0, X] \times [0, T]$. The test function is multiplied by Eq.~\ref{eq:wsindy-general-wtt} and used to perform integration by parts. Assume the following notation for the inner product of real valued functions
    \begin{equation}\label{innp}
        \langle u(x,t), v(x,t) \rangle = \int_0^T \int_0^X u(x,t) v(x,t)\ dx\ dt. 
    \end{equation}

    The inner product of the potential dynamics and test function, i.e. the weak form of Eq.~\ref{eq:wsindy-general-wtt}, is 
    \begin{equation}\label{eq:wsindy-general-wtt22}
        \langle \psi(x,t), w_{,tt}(x,t) \rangle = \langle \psi(x,t), F(w) \rangle. 
    \end{equation}

    On repeating this projection for an ensemble of test functions $\psi^{(k)}(x,t),\ k=1,\cdots, K$, one finds  
    \begin{equation}\label{eq:wsindy-general-wtt44}
        \langle \psi^{(k)}(x,t), w_{,tt}(x,t) \rangle = \langle \psi^{(k)}(x,t), F(w) \rangle. 
    \end{equation}

    After using integration by parts on the left hand side, for each test function, the resulting system is 
    \begin{equation}
        \langle \psi^{(k)}_{,tt}(x,t), w(x,t) \rangle = \langle \psi^{(k)}(x,t), F(w) \rangle, \quad k=1,\cdots, K.
    \end{equation}

    Note that the test functions implemented in the WSINDy for PDEs algorithm smoothly decay to zero at the boundaries of the time and space domains, so the boundary terms from the integration by parts vanish. Additionally, the formulation in \cite{MessengerWeak2021} assumes that $F$ is a linear combination of partial derivatives and nonlinear functions of $w$ chosen by the user \textit{a priori} and known as the library of potential terms. In other words, it is assumed that the PDE can be written in the form 
    \begin{equation} w_{,tt} = F(w) = \sum_{j=1}^J c^{(j)} g^{(j)}(w), \label{eq:wsindy-linearcombo-form}\end{equation}

    where each $g^{(j)}(w),\ j=1, \cdots, J$ is a term in the library of trial operators and $\mathbf{c} = (c^{(1)}, \cdots, c^{(J)})^T$ is the vector of coefficients to solve for. In the present study, the $g^{(j)}(w)$ are $\{w_{,t}, w_{,x}, w_{,xx}, w_{,xxx}, w_{,xxxx}, w, a \}$, where $a$ represents a constant term. The integration by parts is also conducted for the right-hand side of the system such that
    \begin{equation}
        \langle \psi^{(k)}_{,tt}(x,t), w(x,t) \rangle = \bigg\langle \psi^{(k)}_{,*}(x,t), \sum_{j=1}^J c^{(j)} \tilde{g}^{(j)}(w) \bigg\rangle, \quad k=1,\cdots, K. 
    \end{equation}

    where $\psi^{(k)}_{,*}(x,t)$ denotes the partial derivatives and the correct sign on $\psi^{(k)}(x,t)$ necessary to perform the integration by parts for each $g^{(j)}(w)$ term in the library and each $\tilde{g}^{(j)}(w)$ is the library term after performing the integration by parts to move the partial derivatives to the test functions $\psi^{(k)}(x,t)$. 
    
    Then, the system is discretized by evaluating the library terms and left-hand side term on the data grid $(\mathbf{x},\mathbf{t}) \in [0, X]^{N_x}\! \times [0, T]^{N_t}$ affiliated with $\mathbf{W}$. Given the discrete forms $\mathbf{U}$ and $\mathbf{V}$ of the real-valued functions $u(x,t)$ and $v(x,t)$, the discretized inner product may be defined by 
    \begin{equation}
        \langle \mathbf{U}, \mathbf{V} \rangle_d = \frac{X}{N_x} \frac{T}{N_t} \sum_{m=1}^{N_t} \sum_{n=1}^{N_x} U_{nm} V_{nm}, \quad U_{nm} = u(x_n,t_m), \, V_{nm} = v(x_n,t_m),
    \end{equation}
    whereby numerical approximations of the integrals are computed. As such, the discretized system can be written as,
    \begin{equation}
        \langle \mathbf{\Psi}^{(k)}, \mathbf{W} \rangle_d = \bigg\langle \mathbf{\Upsilon}^{(k)}, \sum_{j=1}^J c^{(j)} \tilde{g}^{(j)}(\mathbf{W})\bigg\rangle_d, \quad k=1,\cdots, K, 
    \end{equation}
    
    wherein ${\Psi}_{nm}^{(k)}=\psi^{(k)}_{,tt}(x_n,t_m)$ and ${\Upsilon}_{nm}^{(k)}=\psi^{(k)}_{,*}(x_n,t_m)$. This can also be written in matrix form as
    \begin{equation}
        \mathbf{b} = \mathbf{G} \mathbf{c}, \quad  \mathbf{b} \in \mathbb{R}^{K}, \, \mathbf{c} \in \mathbb{R}^{J},
        \label{eq:wsindy-matrixform}
    \end{equation}
    where
    \begin{equation}
    \begin{cases}
        {b}_{k} = \langle \mathbf{\Psi}^{(k)}, \mathbf{W} \rangle_d \\*[1mm]
        {G}_{k,j} = \bigg\langle \mathbf{\Upsilon}^{(k)}, \tilde{g}^{(j)}(\mathbf{W})\bigg\rangle_d
    \end{cases}\!\!\!\!\!\!\!\!\!, \,\, k=1,\cdots, K,\,\, j=1,\cdots, J,
    \label{eq:wsindy-matrixformdetails}
    \end{equation}
    
    and $\mathbf{b}\in \mathbb{R}^K$ is the left-hand side vector, $\mathbf{G}\in \mathbb{R}^{K\times J}$ is the library matrix, and $\mathbf{c} \in \mathbb{R}^J$ is the vector of coefficients.

    \subsubsection{Solving for $\mathbf{c}$ with Sparse Regression}

    From Eq.~\ref{eq:wsindy-matrixform}, the goal is to find the vector of coefficients $\mathbf{c}$ which identifies the active terms in the library and defines the discovered PDE. This will be accomplished by leveraging the sparse nature of most physical principles. The WSINDy approach employs a modified sequential-thresholding least squares (MSTLS) algorithm to solve for a sparse solution $\mathbf{c}$ \cite{MessengerWeak2021}. An important choice in sparsity promoting methods is the value of the sparsity threshold hyperparameter, which balances promoting sparsity in the coefficient vector as to avoid underfitting with eliminating important terms. As part of the MSTLS algorithm, the optimal sparsity parameter $\hat{\lambda}$ is optimized by minimizing the loss function 
    \begin{equation}
        \mathcal{L}(\lambda) = \frac{\Vert\mathbf{G}(\mathbf{c}_\lambda - \mathbf{c}_{LS})\Vert_2}{\Vert\mathbf{G}\mathbf{c}_{LS}\Vert_2} + \frac{\Vert \mathbf{c}_\lambda \Vert_0}{J}, 
    \end{equation}
        
    where $\mathbf{G}$ is the library matrix described in Eq.~\ref{eq:wsindy-matrixformdetails}, $\mathbf{c}_\lambda$ is the output of an iteration in the MSTLS algorithm, $\mathbf{c}_{LS} = (\mathbf{G}^T \mathbf{G})^{-1} \mathbf{G}^T \mathbf{b}$ is the least squares solution, and $\Vert \mathbf{c}_\lambda \Vert_0$ is the number of non-zero terms in $\mathbf{c}_\lambda$. See \cite{MessengerWeak2021} for additional details. The first term enforces a match between the least squares solution and the coefficient vector, and the second term enforces that the coefficient vector is sparse.

    \subsubsection{Implementation Details}
    
    The implementation in \cite{MessengerWeak2021} has several hyperparameters needed in the algorithm including the choice of library (in terms of the trial differential and nonlinear operators), the search space for the sparsity threshold, the support of one-dimensional test functions, and the query points used to evaluate the test functions as elucidated in the following. The algorithm includes a number of built-in methods to automatically choose several of the hyperparameters and to rescale the data. These are briefly explained in \ref{sec:wsindydetails}. Table~\ref{tab:wsindyparams} also lists the hyperparameter values used in this study which were determined by the relevant sub-algorithms. Specifically, the test functions are piecewise polynomial functions which decay to zero at the boundaries and are separable in time-space. More specifically,
    \begin{equation}\label{reftest}
    \begin{aligned}
    \psi(x,t) &= \phi_1(x)\phi_2(t),\\
        \phi_1(x) &= 
        \begin{cases}
            C_1(x-a)^{p_x}(b-x)^{p_x} \quad a=0 < x < b=X \\
            0 \hspace*{1.25in}\, \text{ otherwise}
        \end{cases}\!\!\!\!\!\!, \quad
        C_1 = \frac{1}{p_x^{2p_x}} \bigg(\frac{2p_x}{X} \bigg)^{2p_x}, \\
        \phi_2(t) &= 
        \begin{cases}
            C_2(t-a)^{p_t}(b-t)^{p_t}\,\, \quad a=0 < x < b=T \\
            0 \hspace*{1.25in} \text{ otherwise}
        \end{cases}\!\!\!\!\!\!,  \quad
        C_2 = \frac{1}{p_t^{2p_t}} \bigg(\frac{2p_t}{T} \bigg)^{2p_t}.  
    \end{aligned}
    \end{equation}
    
    The separability means that the weak form can be computed efficiently by a sequence of one dimensional integrals \cite{MessengerWeak2021}. In the test functions, the polynomial order $p_x\text{ or } p_t$ controls the decay rate towards the domain boundaries. Values for $p_x, p_t,$ and the number of points used to evaluate the support of $\phi_1(x)\text{ and } \phi_2(t)$ are chosen automatically within the analysis based on the spectrum of the data and a corner finding algorithm. These hyperparameters determine how many test functions are used over the space and time domains. There is also an automatic rescaling procedure for the data and spatio-temporal dimensions to regularize the problem and improve the equation discovery for very high or low magnitude data where computing the nonlinear library terms may lead to a large condition number for $\mathbf{G}$ and inhibit accurate discovery of $\mathbf{c}$. The rescaled data is used in the computations of $\mathbf{b}$ and $\mathbf{G}$ and the resulting coefficient vector $\mathbf{c}$ is converted back to the original scaling once identified. It should also be mentioned that the implementation is restricted to the case where each test function $\psi^{(k)}(x,t)$ in Eq.~\ref{eq:wsindy-general-wtt44} is a translation of the reference test function $\psi(x,t)$ in Eq.~\ref{reftest} such that $\psi^{(k)}(x,t) = \psi(x_k-x, t_k-t)$ for a set of query points ${(x_k,t_k)} \in [0, X] \times [0, T]$, $k = 1,2,\ldots,K$. In this setting, the inner products in Eqs.~\ref{eq:wsindy-general-wtt22}-\ref{eq:wsindy-matrixformdetails} take the form of convolution and thus the fast Fourier transform algorithm is leveraged to expedite the integration process.

\subsection{Ensemble Weak SINDy for PDEs}
\label{sec:ensemble-wsindy}

   This work also implements an ensemble version of the WSINDy for PDEs algorithm to improve the equation discovery in the presence of measurement noise and gain information about the uncertainty in the coefficients and terms in the discovered PDEs stemming from noise and variations in the data. Since the WSINDy algorithm assumes the data is on a uniform grid, subsamples of the full data set are taken with a fixed step size between data points. The user defines a maximum downampling rate and for each one, the starting point of the data subsets is looped through to use all the available data. Thus, the subsampling takes place on the data before being passed to the WSINDy algorithm to perform the equation discovery. This is outlined in Algorithm~\ref{alg:ensemble-wsindy}. In this study, the ensembling of the full data set only takes place in time. Since for each data set the number of time points is by an order of magnitude larger than the number of spatial points. Thus by sampling in time a larger number of subsets can be taken while still maintaining enough data for accurate equation discovery.

    Note, other work \cite{FaselEnsemble2022} has investigated alternate ensemble approaches with SINDy methods by (1) bootstrapping the data and (2) bagging the terms in the library. Rather than using data bootstraps, our approach chooses the data subsets to preserve the uniform grid requirement of the data. In this study, ensembling on the library terms is not performed since the library starts with a small number of terms and taking subsets does not result in reliable equation discovery. However, the library bagging could be helpful in settings that start from a very large library to increase the efficiency of the equation discovery procedure.

    \begin{algorithm}
    \caption{Data Ensemble WSINDy for PDEs}\label{alg:ensemble-wsindy}
     \hspace*{\algorithmicindent} \textbf{Input}: $\mathbf{W}$, maximum downsampling rate, WSINDy \\
     \hspace*{\algorithmicindent} \textbf{Output}: $\mathbf{c}$ per ensemble
        \begin{algorithmic}
        \For {$d = 1:max\_ds\_rate$}
            \For {$i = 1:dsi$}
                \State $\tilde{\mathbf{W}} = W_{nm}$ $\quad \forall n,\ m = \{i, i+d, i+2d, i+3d,\cdots  \le N_t \}$ 
                \State $\mathbf{c} = $ WSINDy($\tilde{\mathbf{W}}$)
                \State Return $\mathbf{c}$
            \EndFor
        \EndFor
        \end{algorithmic}
    \end{algorithm}

\section{Results}
\label{sec:results}

\subsection{Equation Discovery}
\label{sec:eqdiscovery}

    The WSINDy algorithm was run with the two experimental data windows shown in Fig.~\ref{fig:exp-data} (b) and a library containing time derivatives up to order two, spatial derivatives up to order four, linear terms in the field variable $w$, and a constant term ($a$), that is library $ = \{w_{,t}, w_{,x}, w_{,xx}, w_{,xxx}, w_{,xxxx}, w, a \}$. The second time derivative $w_{,tt}$ is set as the left hand side of the equation because this term is associated with inertia and often present in equations governing the dynamics of beam-like structures. 

\begin{figure}[h]
        \centering
        \includegraphics[width=0.55\textwidth]{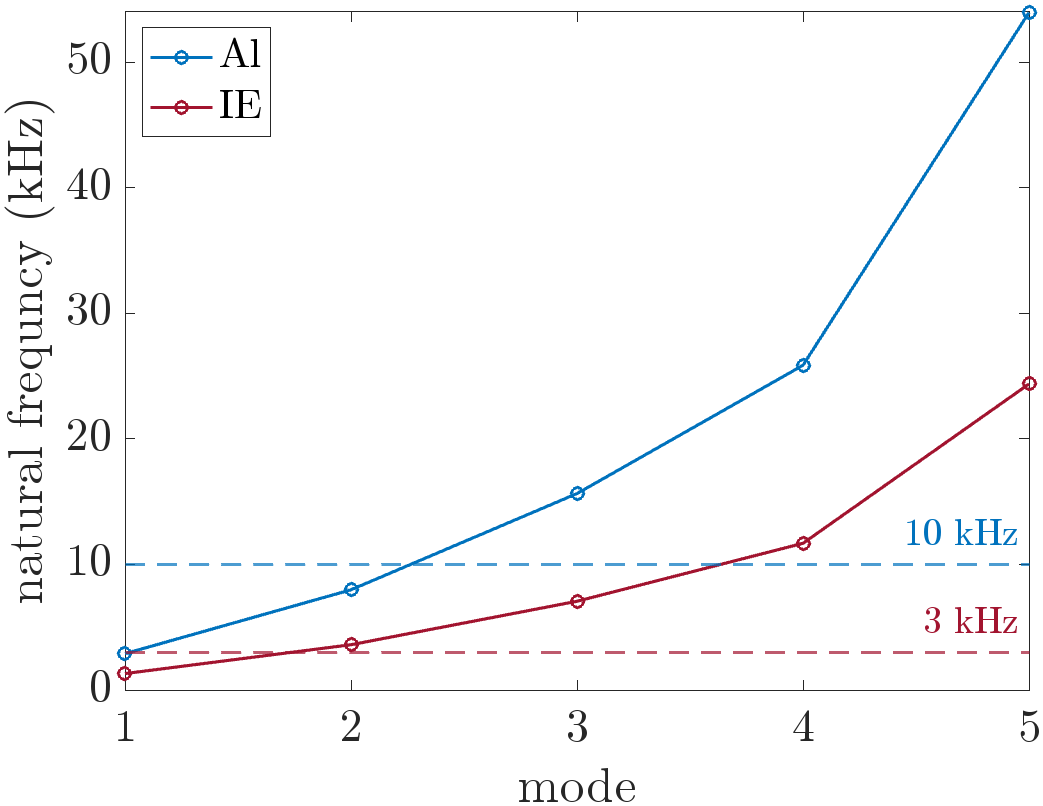}
        \caption[Natural frequency vs mode]{Natural frequency vs mode for Al and IE beams computed using the testing configuration, the geometry of specimens and $E$ values in Table~\ref{tab:eqdiscovery}. The Al results are shown in blue and the IE results in red. The input frequencies are marked with the dashed lines and fall within the first one to two modes for each specimen. }
        \label{fig:modes}
    \end{figure}

    The partial differential equations (PDEs) discovered for both the Al and IE data were of the form 
    \begin{equation}
        w_{,tt} = - \alpha w_{,xxxx}.
    \end{equation}
    Table~\ref{tab:eqdiscovery} summarizes the discovered equations. In other words, of the terms in the library, the WSINDy algorithm identified only the $ w_{,xxxx}$ term as being active in the data. This model form can be interpreted as the Euler-Bernoulli beam model for a homogeneous beam with no external loading. This beam model describes the deflection of a beam in terms of displacement (in the homogeneous case the deflection can also be described in terms of velocity). The Euler-Bernoulli model form is plausible for these tests since the specimens are slender (cross section is more than 10 times smaller than the length) and the input signals are in the low frequency/long wavelength regime as they fall within the first few modes for the geometry and materials of the respective specimens (Fig.~\ref{fig:modes}). It is important to notice that the approach did not ``fit" an Euler-Bernoulli model to the data, rather the data, along with WSINDy, led to the Euler-Bernoulli model as the most plausible. The WSINDy relative residuals, computed according to $\Vert \mathbf{b}-\mathbf{Gc}\Vert_2 / \Vert \mathbf{b}\Vert_2$, for the discovered PDEs are $17.1$\% and $21.5$\% for the Al and IE, respectively. This indicates that to the leading order, the dynamics of the data are captured by the Euler-Bernoulli beam model, but some aspects of the behavior are not described by this model.  

    \begin{table}[h!]
        \centering
        \begin{tabular}{|l|l|l|l|l|}
        \hline
        & \textbf{Discovered PDE} & \textbf{WSINDy Relative Residual} & \textbf{Estimated E (Pa)} & \textbf{E \% Error} \\ \hline
        \hline
         Al & $w_{,tt} = -58.5218w_{,xxxx}$ & $0.171$ & $6.3206e+10$ & $8.3967\%$ \\ \hline
         EI & $w_{,tt} = -0.497308 w_{,xxxx}$ & $0.215$ & $9.6292e+08$ & N/A \\ \hline
        \end{tabular}
        \caption[Equation discovery results.]{Equations discovered using WSINDy. The relative residual is reported by WSINDy and is computed according to $\Vert \mathbf{b}-\mathbf{Gc}\Vert_2 / \Vert \mathbf{b}\Vert_2$. The estimated $E$ is computed from the coefficient in the discovered PDE and geometry of the sample. For aluminum, the percent error is computed with respect to the nominal Young's modulus \cite{ASMHandbookAl}.}
         \label{tab:eqdiscovery}
    \end{table}
    
    In the Euler-Bernoulli interpretation of the discovered model, the coefficient on the $w_{,xxxx}$ term is $\alpha = EI/\rho A$, where $E$ is the Young's modulus, $A$ is the cross section area, $I$ is the second moment of area, and $\rho$ is the mass density. Since $A, I, \rho$ are known from the geometry and mass of the specimen, the coefficient and known values can be used to also recover an estimate for $E$, which is not measured directly from the specimen. For the Al specimen, the recovered value for $E$ is $6.3206e+10$ \unit{\pascal}, which is within $10\%$ of the nominal value for aluminum of $6.9e+10$ \unit{\pascal} \cite{ASMHandbookAl}. The estimate of $E$ for the IE specimen is $9.6292e+08$ \unit{\pascal}. While there is not an established nominal value for the Young's modulus of the composite, this estimate is within the bounds of Young's modulus values for the individual IDOX and Estane components. Young's modulus of the Estane binder is about $2.3 e+06$ \unit{\pascal} at room temperature \cite{Estane} and Young's modulus of the IDOX crystals has been measured at $2.33 \pm 0.22 e+10$ \unit{\pascal} via nanoindentation \cite{YeagerDevelopment2019, BurchNanoindentation2017}. Previous studies which tested the IE composite under quasi-static loading conditions at ambient temperature, for specimens with different manufacturing process, have reported that the Young's modulus for the composite is between $3.77e+06$ \unit{\pascal} and $6.273 e+08$ \unit{\pascal} \cite{HermanComposite2021, LiuPBX95012020, YeagerDevelopment2019}. The Young's modulus recovered for the present specimen is somewhat larger than these previously reported values for the composite, but still consistent with the previous results given differences in the manufacturing and testing procedures. Thus, for both specimens the discovered PDE and associated Young's modulus are physically plausible. 

\subsection{Comparisons with Finite Element}
\label{sec:validation}

    In order to compare the equations discovered with WSINDy to another method, a custom elastodynamic finite element code was used to solve the beam equation PDE and compare the results to the experimental velocity data. The finite element code follows a standard finite element formulation \cite{HughesFinite2000} for beam elements and Hermite cubic test functions, and is implemented in MATLAB \cite{MATLAB}. Briefly, the finite element code simulates a beam with prescribed displacement and rotation boundary conditions at $x=0$ and $x=L$, by solving 
    \begin{equation}
        \mathbf{M} \mathbf{\Ddot{d}} + \mathbf{K} \mathbf{d} = \mathbf{F},
        \label{eq:fem}
    \end{equation}
    in time domain, where $\mathbf{M}$ is the global mass matrix, $\mathbf{K}$ is the global stiffness matrix, $\mathbf{d}$ contains the nodal displacements and rotations excluding the boundary nodes, $\mathbf{\Ddot{d}}$ has the nodal acceleration terms also excluding the boundary nodes, and $\mathbf{F}$ is the global force vector due to the prescribed displacement boundary conditions at the boundary nodes. 
    
    To model the experimental configuration as closely as possible, the spatial domain and geometry for the simulations were set to match the experimental specimens. All elements are the same length and match the spacing of the sampling points in the experimental setup; $dx = 0.5$ \unit{\milli\meter}. Thus, $194$ beam elements were used for the Al simulation and $112$ were used for the IE. The recovered Young's modulus is treated as an input for the simulations. Since the experimental data is collected at short intervals, the time step in the simulation was set to match the sampling rates of the experimental data; $dt = 1.600e-07$ \unit{\second} for Al and $dt = 3.195e-07$ \unit{\second} for IE. Each simulation starts at $t=0$ where the initial condition $w(x, 0) = w_{,t}(x, 0) = 0$ is satisfied, but the results are presented only for the window matching the experimental data used for equation discovery since it is the time span where the PDE is expected to hold. To handle the temporal component, the finite element equation, Eq.~\ref{eq:fem} is integrated in time using the Newmark's method with the trapezoidal rule. The prescribed displacements at the boundaries are enforced with the experimental measurements at the first and last sampling point in each data set. The prescribed rotation at each boundary is computed from the experimental measurements as follows, the first 25 points and last 25 points in space are fit with a third order Fourier series. The Fourier series is used to analytically compute the spatial derivative $w_{,x}$ at $x=0$ and $x=L$, which is applied as the prescribed rotation at the boundaries. This is repeated for each time point. These two steps result in a time series of prescribed displacement and rotation values for the two boundary locations. Then, to compute the acceleration of the prescribed boundary values, the second derivative of the time series is computed numerically with a differencing scheme. At each Newmark iteration the prescribed boundary values and force vector are updated. Note, the spatial and temporal grids for the simulations were chosen to avoid interpolation in the comparison with the experimental data.
    
\begin{figure}[h]
        \centering
        \includegraphics[width=0.99\textwidth]{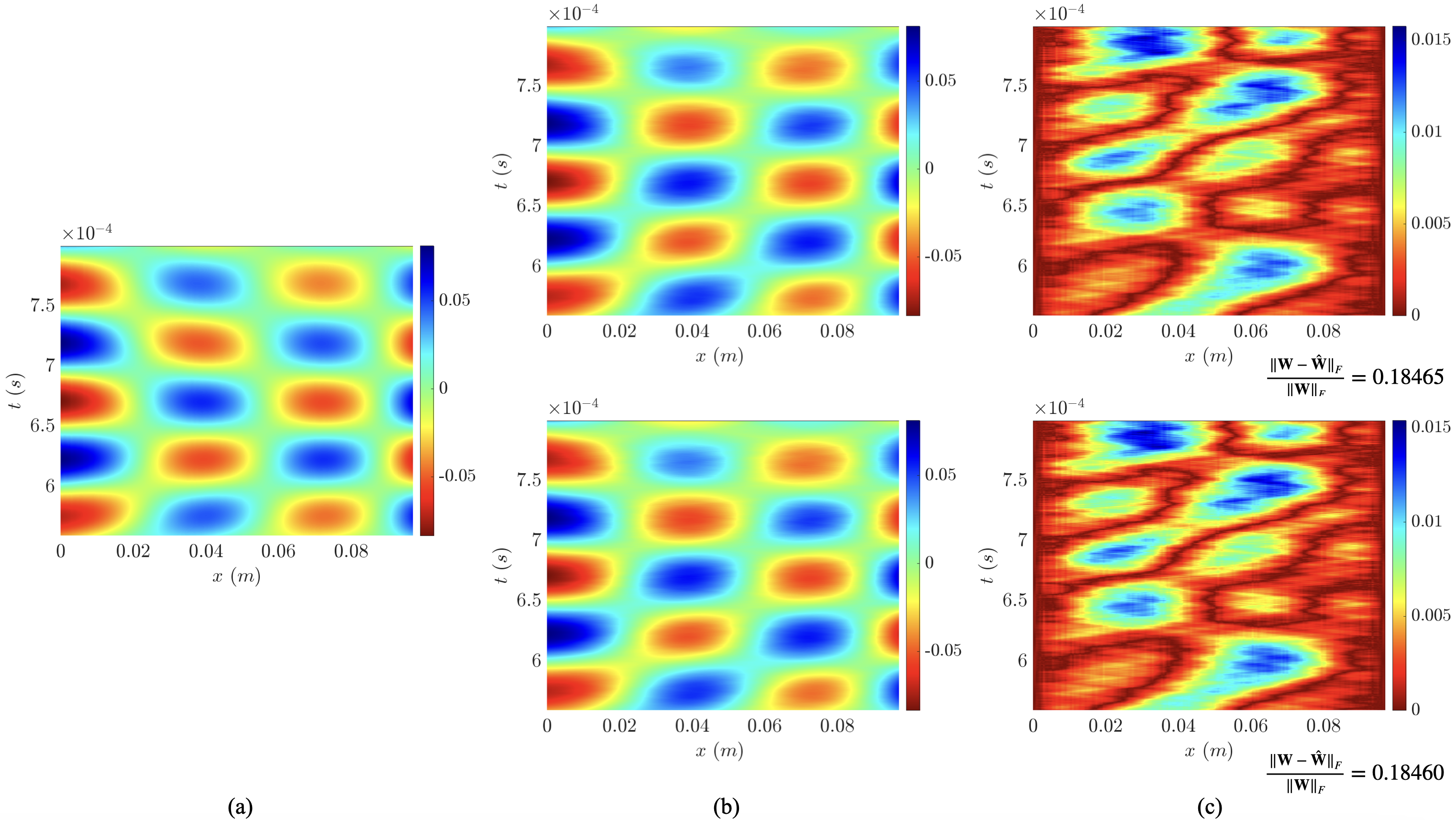}
        \caption[Al FEM]{Al FEM results for two values of $E$. The experimental field is shown in (a), field simulated with FEM is shown in (b), and a error field computed via $|\mathbf{W}-\hat{\mathbf{W}}|$ is shown in (c). The top row corresponds to $E=6.3206e+10$ Pa found from the WSINDy coefficient and the bottom row corresponds to $E=6.33056e+10$ Pa found from minimizing the relative Frobenius error.}
        \label{fig:al-fem}
    \end{figure}

\begin{figure}[h]
        \centering
        \includegraphics[width=0.99\textwidth]{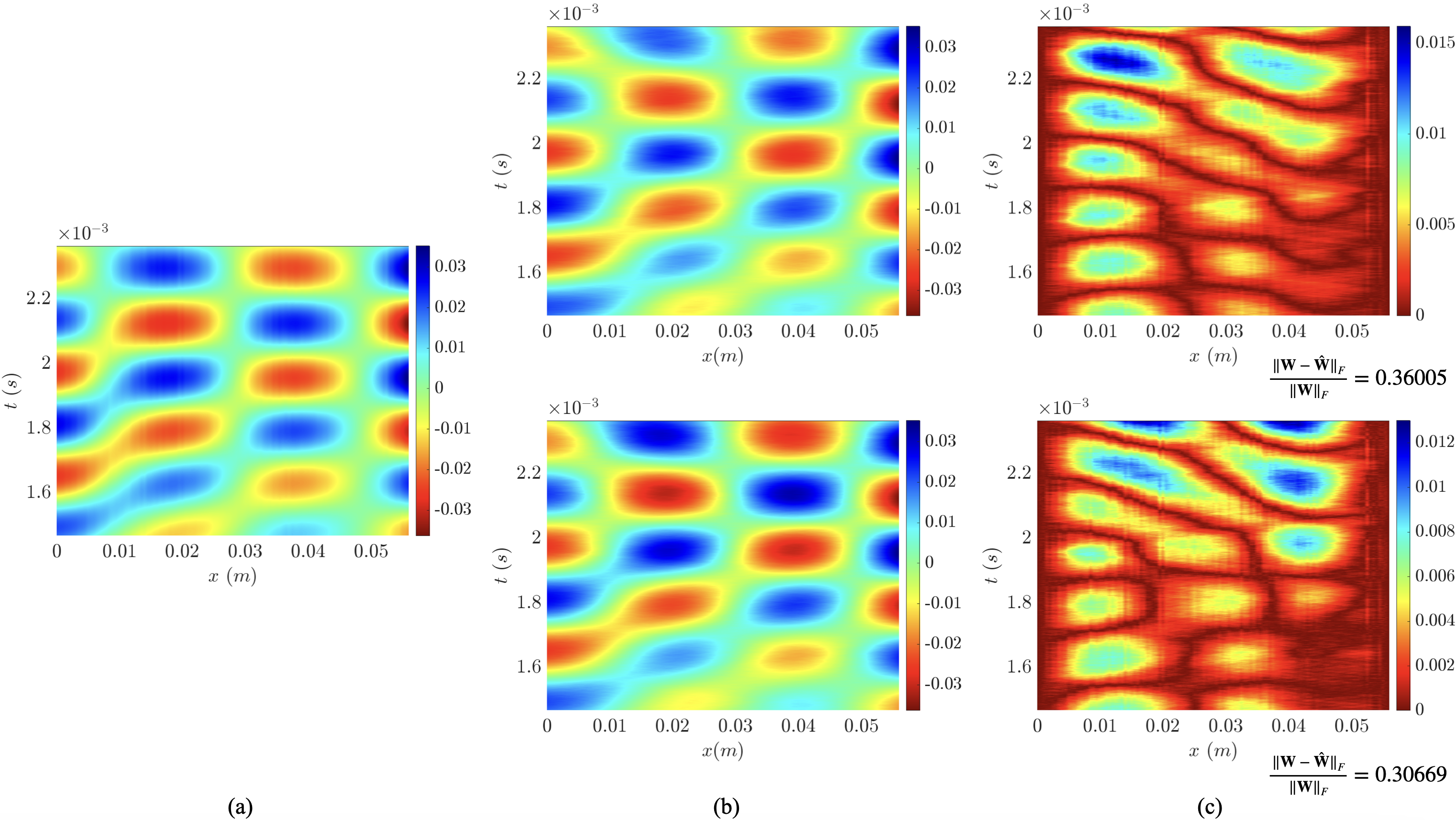}
        \caption[IDOX FEM]{IE FEM results for two values of $E$. The experimental field is shown in (a), field simulated with FEM is shown in (b), and a error field computed via $|\mathbf{W}-\hat{\mathbf{W}}|$ is shown in (c). The top row corresponds to $E=9.6292e+8$ \unit{\pascal} found from the WSINDy coefficient and the bottom row corresponds to $E=1.05912e+9$ Pa found from minimizing the relative Frobenius error.}
        \label{fig:idox-fem}
    \end{figure}

    The finite element code was first used to simulate the equations discovered in Section~\ref{sec:eqdiscovery} with their respectively recovered $E$ values. The top row of Figs.~\ref{fig:al-fem} and ~\ref{fig:idox-fem} show the simulated full-field data and point-wise error field computed according to $|\mathbf{W}-\hat{\mathbf{W}}|$, where $\mathbf{W}$ is the experimental data field as before and $\hat{\mathbf{W}}$ is the velocity field from the FEM simulation. Also for each simulation a relative Frobenius error is computed to give a scalar summary of the error field, according to $\Vert \mathbf{W} - \hat{\mathbf{W}} \Vert_F / \Vert \mathbf{W} \Vert_F$, where $\Vert \cdot \Vert_F$ is the Frobenius matrix norm. The relative Frobenius error for the Al simulation is $0.18465$ and $0.36005$ for the IE simulation. These errors are sightly higher than the WSINDy relative residuals, which is likely due the rotation and acceleration boundary information being computed from fits of the experimental data at the boundaries.

\begin{figure}[h!]
        \centering
        \subfloat[]{\includegraphics[width=0.46\textwidth]{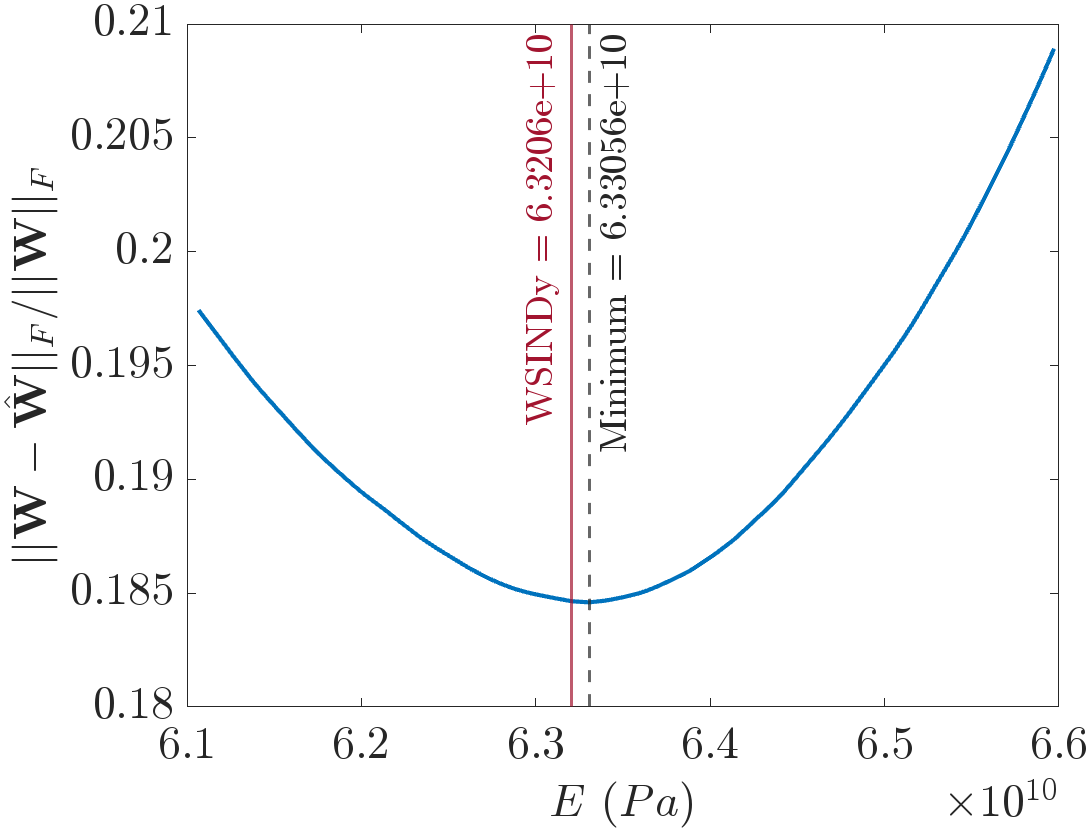}}
        \hfill
        \subfloat[]{\includegraphics[width=0.46\textwidth]{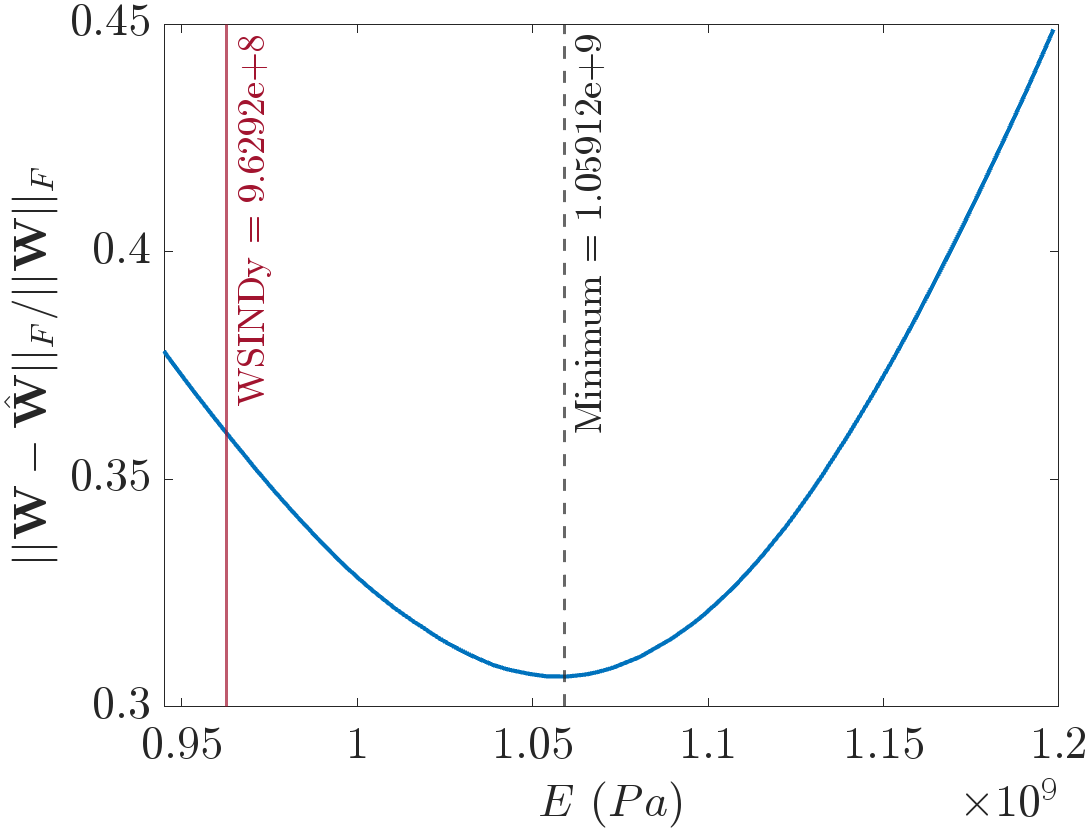}}
        \caption[Manual optimization of E]{Minimization of the relative Frobenius error vs $E$ for Al (a) and IE (b). The $E$ value which minimizes the error is marked with a black dashed line and $E$ recovered from the WSINDy coefficient is shown with a solid red line.}
        \label{fig:manualopt}
    \end{figure}

    Next, the finite element code was used to investigate if the relative Frobenius error between the experimental and simulated fields could be improved upon by assuming the Euler-Bernoulli model form discovered by WSINDy and varying $E$ to minimize the error. A grid of $E$ values around the solution from the WSINDy coefficient was chosen and used to run a set of simulations and compute the errors. For the Al data, the simulation was run for $150$ values of $E$ sampled from $[6.1 e+10,\  6.6 e+10]$ \unit{\pascal} and was run for $150$ values sampled from $[9.2e+08,\ 1.2 e+09]$ \unit{\pascal} for IE. Fig.~\ref{fig:manualopt} shows plots of the error between the simulations and the experimental data against $E$. The error is minimized to $0.18460$ at $6.33056e+10$ \unit{\pascal} and $0.30669$ at $1.05912e+09$ \unit{\pascal} for the Al and IE data, respectively. The simulation results of from the finite element code with these optimized values of $E$ are shown in the bottom row of Figs.~\ref{fig:al-fem} and ~\ref{fig:idox-fem}. In the Al case, the WSINDy and optimized values of $E$ are very close and the impact on the simulated data fields and error fields is negligible. However, in the IE case, the $E$ values from the two methods differ by $9.62e+07$ \unit{\pascal} and the optimized value decreases the error by about $5.3\%$ compared to the WSINDy value. By visual inspection, the IE velocity fields in Fig.~\ref{fig:idox-fem} (b) are similar. Interestingly though, the error fields have the highest levels of error in different locations for the two values of $E$, as can be seen in Fig.~\ref{fig:idox-fem} (c). For the WSINDy $E$ value, the highest error is concentrated in the late time, early in space quadrant whereas the optimized $E$ value leads to the highest error concentrated in the later times across the spatial domain. Interestingly, the $E$ which minimizes the error is higher than any previous experimental measurements of the Young's modulus for the composite, which calls into question whether the modulus is a physically reasonable quantity. Previous work with similar experimental data on 2D plates has shown that inverting the data based on a model which does not capture all the relevant physics or where the data are just outside the intended regimes of the model, the resulting optimal values can be non-physical to account for the discrepancy in the model and the data ~\cite{XuDeep2023}. To determine if this is the case here or if it is a valid value for $E$ which has not yet been observed in other work, additional follow up testing should be conducted. 

    Overall, the WSINDy and optimized $E$ values are similar and result in similar simulated fields with only a few percentage points difference in the error to the experimental data for both materials. The benefit of using WSINDy here is in the ability to discover the model form and coefficients simultaneously as opposed to assuming a model form and calibrating the material property of interest. Additionally, it is important to note that, based on the design of the WSINDy test functions decaying to zero at the domain boundaries, no information about the boundary conditions is needed for the equation discovery. Whereas for the finite element simulations, significant effort is needed to generate the boundary condition data from experimental measurements. The WSINDy algorithm is also significantly more efficient than running the finite element simulation, especially $150$ times. For the Al and IE results presented in Table~\ref{tab:eqdiscovery}, each equation discovery took about $0.79$ \unit{\second} whereas each finite element simulation takes about $5.6$ \unit{\min} and each set of $150$ finite element simulations approximately $14$ \unit{\hour} to run on a MacBook Pro with an M2 Pro chip. This is not intended as a rigorous comparison of the computational efficiency of the two approaches since optimizing the finite element implementation for speed was outside the scope of this study, but rather to highlight that the two methods reached similar results on different time scales. In settings with higher dimensions, more complex geometry, non-linearities, etc, even with a state of the art finite element code, these timing differences are likely to be even more pronounced. Additionally, for applications where an equation is needed with quick turn around, an equation discovery method such as WSINDy could give an answer faster than performing model selection and calibrating values in a forward model.

\subsection{Ensemble Equation Discovery}
\label{sec:ensemble-eqdiscovery}

\begin{figure}[h]
        \centering
        \subfloat[]{\includegraphics[width=0.32\textwidth]{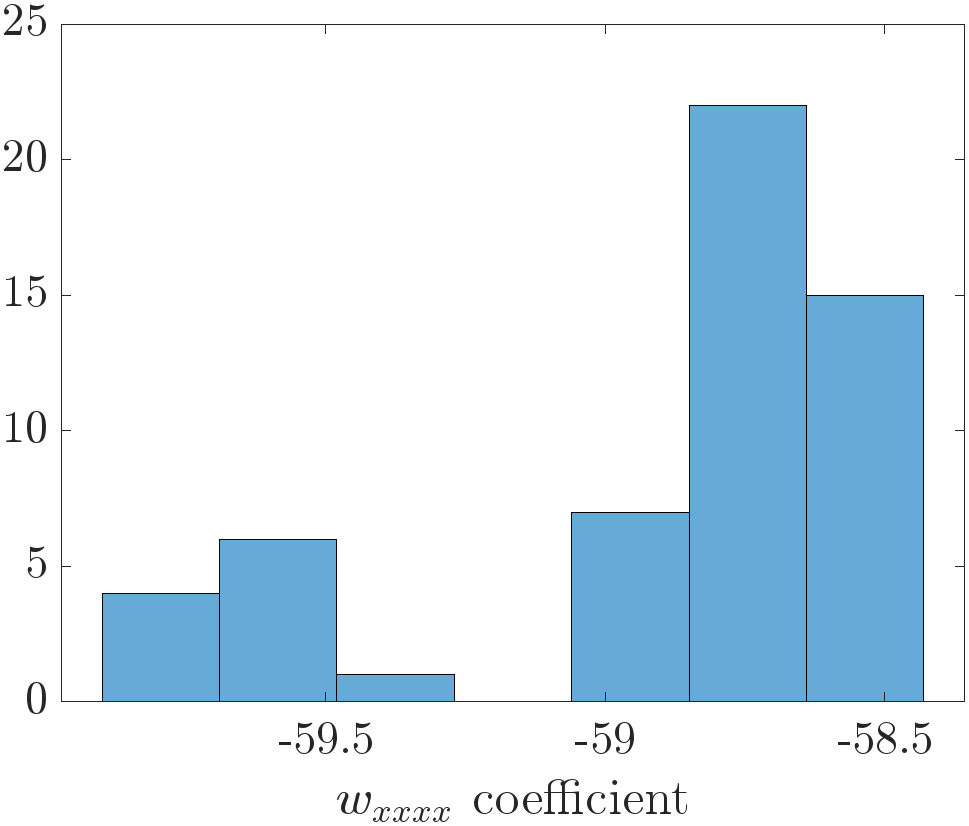}}
        \hfill
        \subfloat[]{\includegraphics[width=0.32\textwidth]{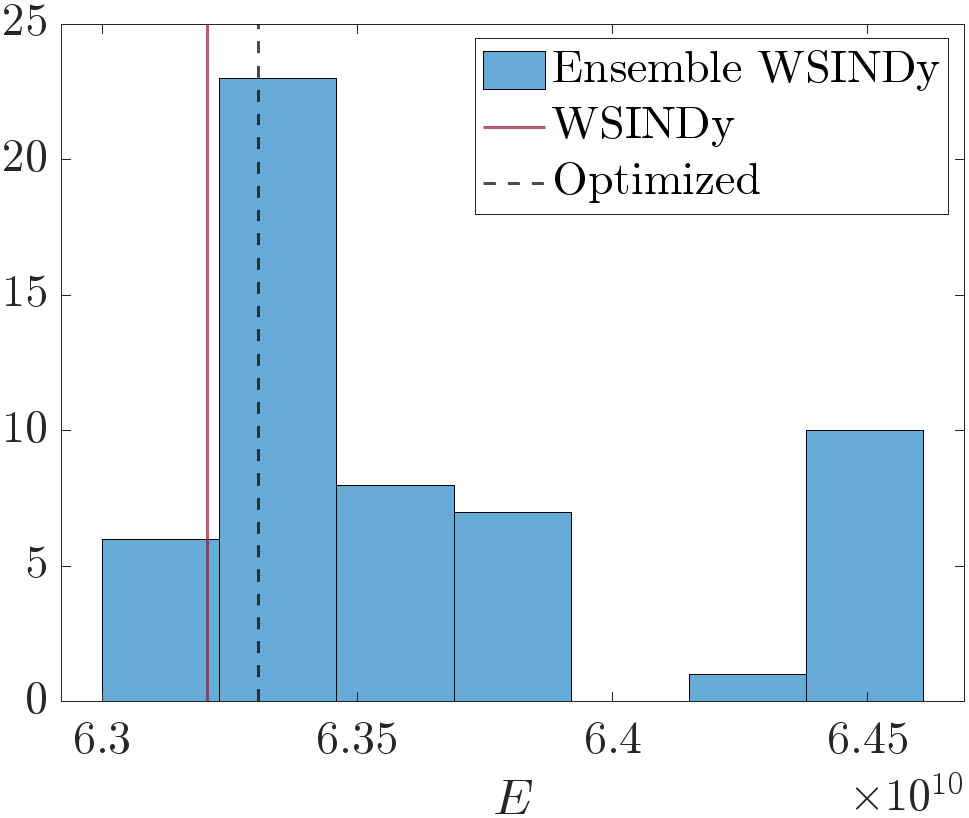}}
        \hfill
        \subfloat[]{\includegraphics[width=0.32\textwidth]{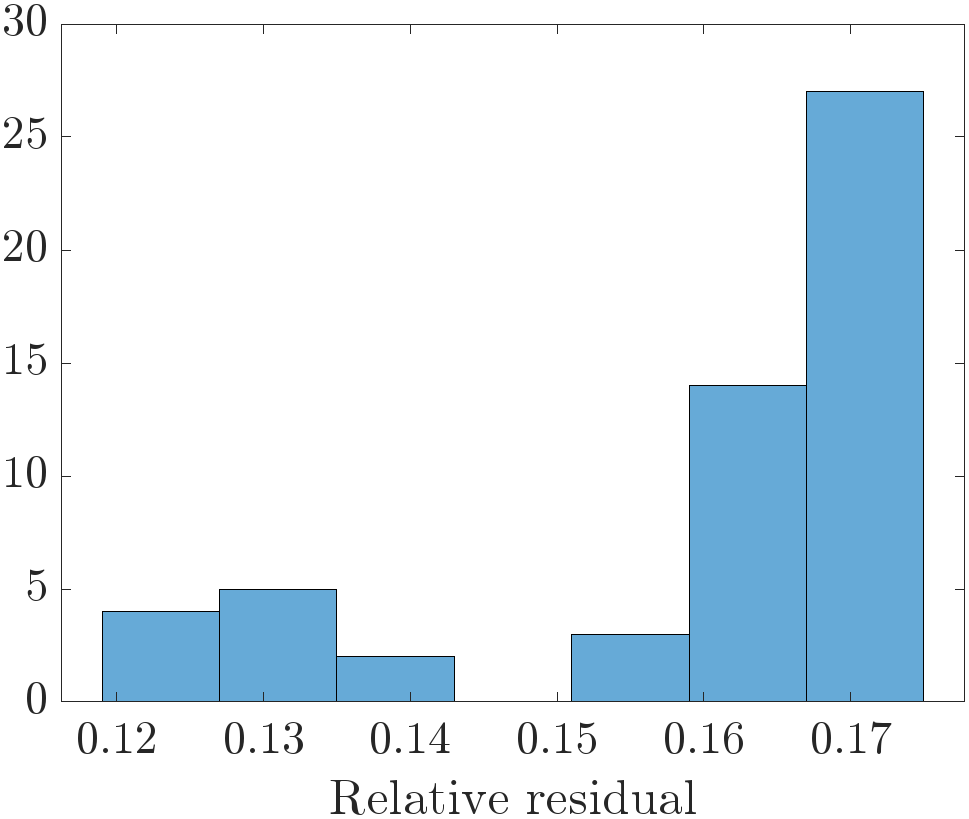}}
        \caption[Al ensemble equation discovery results.]{Histograms of ensemble equation discovery coefficients (a), $E$ values estimated from the coefficients (b), and WSINDy relative residuals (c) for Al experimental data.}
        \label{fig:hist-al}
    \end{figure}

\begin{figure}[h]
        \centering
        \subfloat[]{\includegraphics[width=0.32\textwidth]{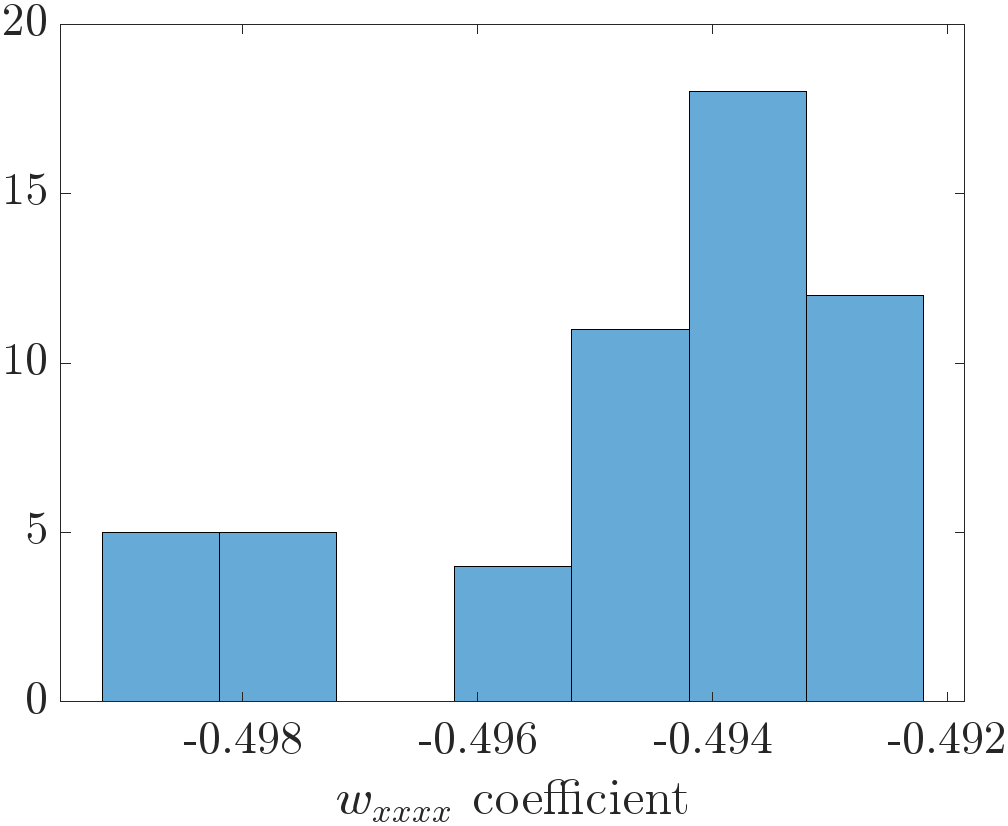}}
        \hfill
        \subfloat[]{\includegraphics[width=0.32\textwidth]{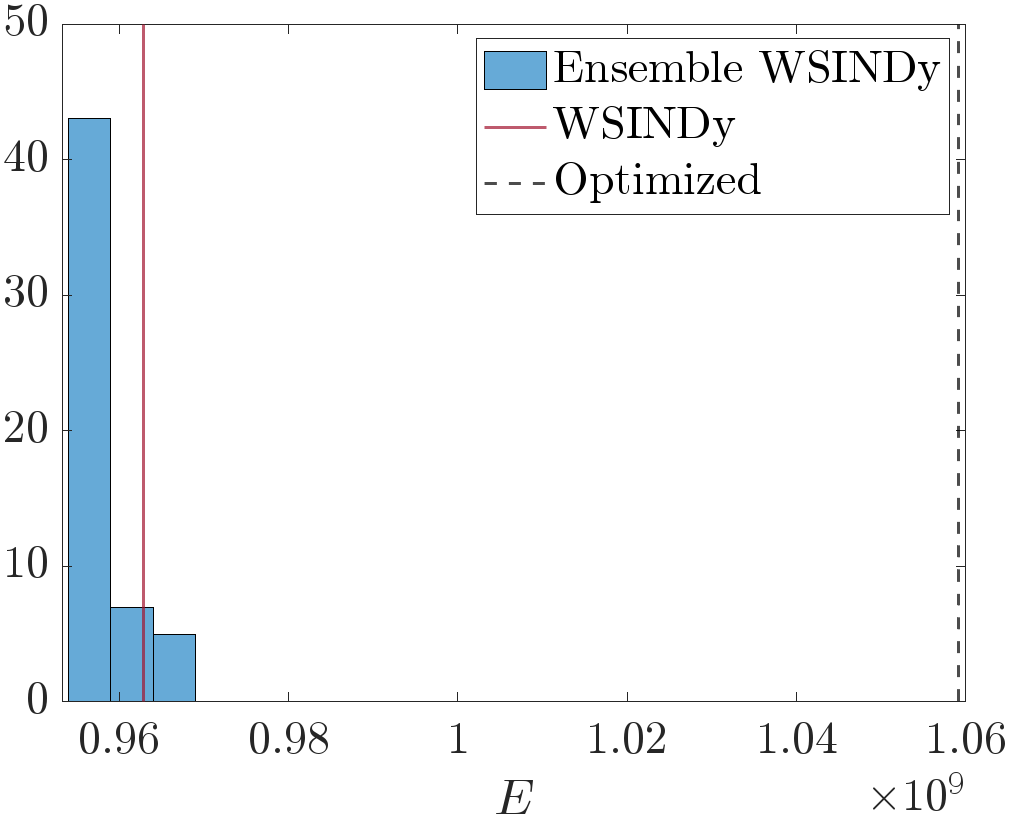}}
        \hfill
        \subfloat[]{\includegraphics[width=0.32\textwidth]{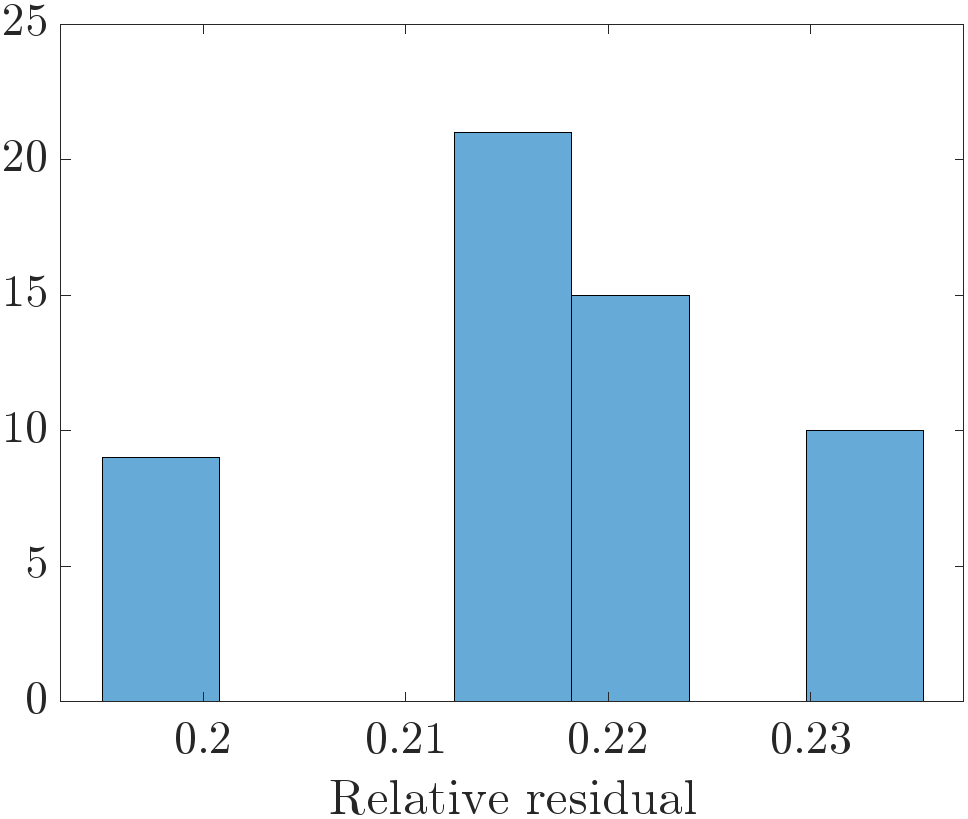}}
        \caption[IE ensemble equation discovery results.]{Histograms of ensemble equation discovery coefficients (a), $E$ values estimated from the coefficients (b), and WSINDy relative residuals (c) for IE experimental data.}
        \label{fig:hist-idox}
    \end{figure}

    The ensemble WSINDy algorithm outlined in Section~\ref{sec:ensemble-wsindy} was run with the same experimental data windows and library as in Section~\ref{sec:eqdiscovery}. A maximum downsampling step of $10$ was used, which resulted in $55$ ensemble datasets and recovered equations for each experiment. For both datasets, all $55$ of the discovered PDEs were of the Euler-Bernoulli form as before. Thus, even taking subsets of the data, the Euler-Bernoulli model captures the dynamics of the system, with some variation in the coefficient on the $w_{,xxxx}$ term. Table~\ref{tab:ensemble-eqdiscovery} and Figs.~\ref{fig:hist-al} and ~\ref{fig:hist-idox} summarize these results.  

    \begin{table}[h]
        \centering
        \begin{tabular}{|l|l|l|l|l|}
        \hline
        & \textbf{Mean E (Pa)} & \textbf{Median E (Pa)} & \textbf{Standard deviation (Pa)} & \textbf{E SMAPE} \\ \hline
        \hline
        Al & $6.3624e+10$  & $6.3424e+10$ & $4.4800e+08$ & $8.1102\%$ \\ \hline
        EI & $9.5769e+08$ & $9.5633e+08$ & $3.5637e+06$ & N/A \\ \hline
        \end{tabular}
        \caption[Ensemble equation discovery results.]{Ensemble equation discovery results. The mean, median, and standard deviation are presented for $55$ values of $E$ which are computed from the coefficient in the ensemble discovered PDEs using the geometry of the sample. For aluminum, the symmetric mean absolute percent error (SMAPE) is computed between the $55$ $E$ values and the nominal Young's modulus and indicates that all values are within $10\%$ of the nominal. }
        \label{tab:ensemble-eqdiscovery}
    \end{table}
    
    As before, estimates for $E$ are recovered using the model coefficient and known values of $A, I, \rho$ for each specimen. The average recovered $E$ values are $6.3624e+10$ and $9.5769e+08$ \unit{\pascal} for the Al and IE tests, respectively. These values are consistent with those found using the WSINDy algorithm applied to the whole dataset (Table~\ref{tab:eqdiscovery}) and fall within the ranges of values from the ensemble equation discovery for both materials. This is illustrated in Figs.~\ref{fig:hist-al} (a) and ~\ref{fig:hist-idox} (a). For Al, the value of the coefficient on the $w_{,xxxx}$ term is in the range $[-59.807, -58.478]$, leading to a range of $[6.3159e+10, 6.4595e+10]$ \unit{\pascal} for the recovered $E$ values. Also, all the ensemble $E$ values are within $10\%$ of the nominal for Al as indicated by the symmetric mean absolute percent error (SMAPE) given in Table~\ref{tab:ensemble-eqdiscovery}. For IE, the coefficient values are in the range $[-0.49894, -0.49286]$ and the recovered $E$ values are in $[9.5430e+08, 9.6607e+08]$ \unit{\pascal}. Thus, taking subsets of the data does not bias the PDE coefficients and resulting $E$ values to be larger or smaller than those discovered with the whole datasets. In the Al case, the $E$ value obtained from minimizing the relative Frobenius error between the experimental and simulated data is also within the range of the ensemble equation discovery results. However, this is not the case for IE where the $E$ found via minimization is larger than the $E$ values found through the WSINDy approach. Lastly, the standard deviation of the ensemble $E$ values for Al is $4.4800e+08$ and is $3.56e37e+06$ \unit{\pascal} for IE. The standard deviations and distributions shown in Fig.~\ref{fig:hist-al} (b) and Fig.~\ref{fig:hist-idox} (b) give a description of the uncertainty in the estimates of $E$ for the two materials. 

    The ensemble WSINDy relative residuals for Al are between $0.123$ and $0.175$ and between $0.195$ and $0.235$ for IE. As for the PDE coefficient and $E$, the relative residuals of the equations discovered over the whole data sets fall within these ranges. This indicates that in performing the ensembling, using smaller subsets of the data does not bias the discovered equations to be more or less accurate than using the entire data set.

\section{Conclusions}
\label{sec:conclusion}

    This work utilizes the WSINDy for PDEs algorithm with experimental data to successfully learn interpretable governing equations. The experimental data are full-field waveforms generated from low frequency laser vibrometry tests on beam-like specimens of aluminum and an IDOX/Estane mock HE composite. The discovered equations are of the recognizable Euler-Bernoulli form which allows the coefficients to be interpreted as plausible estimates for the material Young's modulus. An ensemble version of the WSINDy algorithm provided additional uncertainty information about the resulting Young's modulus values accounting for noise and variation in the data. Additionally, since the recovered model form is simple it can be solved with a finite element approach for the same spatio-temporal grid as the experimental data and against the experimental results. Two benefits of using WSINDy to discover the model over calibrating $E$ in the Euler-Bernoulli model is the ability to learn the model form and coefficients simultaneously, and to avoid generating boundary condition data from noisy experimental measurements. This work also provides an example of analyzing full-field data in the time domain instead of converting to the frequency domain. Overall, this study demonstrates the power of learning unknown dynamics and recovering material properties from experimental data in a data driven framework.

    Since the discovered equations and estimates for Young's modulus are based on data from tests conducted in the low frequency/long wavelength regime, an interesting avenue of future work could be to conduct further tests using higher frequency inputs and perform the equation discovery. For higher frequencies, it is unlikely that the Euler-Bernoulli model would hold, so investigating which model forms would be the most plausible and potentially how the coefficients could be interpreted to gain insights about the materials tested, would be valuable. For additional frequencies close to the long wavelength regime investigated here, i.e. still within the first few modes, the Euler-Bernoulli model form is likely to be the equation discovered by WSINDy. Thus, it could be possible to estimate frequency dependent Young's modulus from a series of low frequency data sets.

\section*{Acknowledgements}

    This work was supported by the Department of Energy, National Nuclear Security Administration, Predictive Science Academic Alliance Program (PSAAP), Award Number DE-NA0003962. The authors would also like to thank Daniel Messenger and David Bortz for helpful discussions regarding the WSINDy algorithm, Jian Song, Conglin Wang, Emi Szabo for assistance conducting the experiments, and Summer Carmelo for manufacturing the IDOX/Estane sample used in the experiments.

\appendix
    
\include{appendix-wsindydetails}


\bibliography{ref}

\end{document}

%% file: appendix-wsindydetails.tex
\section{WSINDy Hyperparameters}
\label{sec:wsindydetails}

    The WSINDy for PDEs algorithm has several hyperparameters, most of which have methods to choose them automatically  \cite{MessengerWeak2021}. This appendix lists the choices and values for the hyperparameters, with short descriptions for clarity. Please see the original publication of the algorithm \cite{MessengerWeak2021} for a more comprehensive description of the hyperparameters and associated methods to choose them. For all of the equation discovery analyses in this study, the piecewise polynomial test functions with automatic selection of the number of points used to define the width of the support for $\phi_1(x)$ and $\phi_2(t)$, $m_x, m_t$ and orders $p_x, p_t$ were used. The algorithm also uses query points ($s_x, s_t$) in the evaluation of the test functions and construction of the linear system $\mathbf{b}=\mathbf{Gc}$ as a discretized convolutional weak form. The default formula for the selection of the query points was used. Additionally, the automatic rescaling procedure and $\tau = 10 e-10$ were used. Finally, for all analyses, the MSTLS algorithm for sparse regression and optimization of the sparsity threshold $\lambda$ was used with a search space of $100$ values between $[1e-10, 1]$ for $\lambda$. The MSTLS algorithm presented in \cite{MessengerWeak2021} enforces a sparse solution for $\mathbf{c}$ while also choosing the optimal sparsity threshold $\hat{\lambda}$.
    
\begin{figure}[h!]
        \centering
        \subfloat[]{\includegraphics[width=0.47\textwidth]{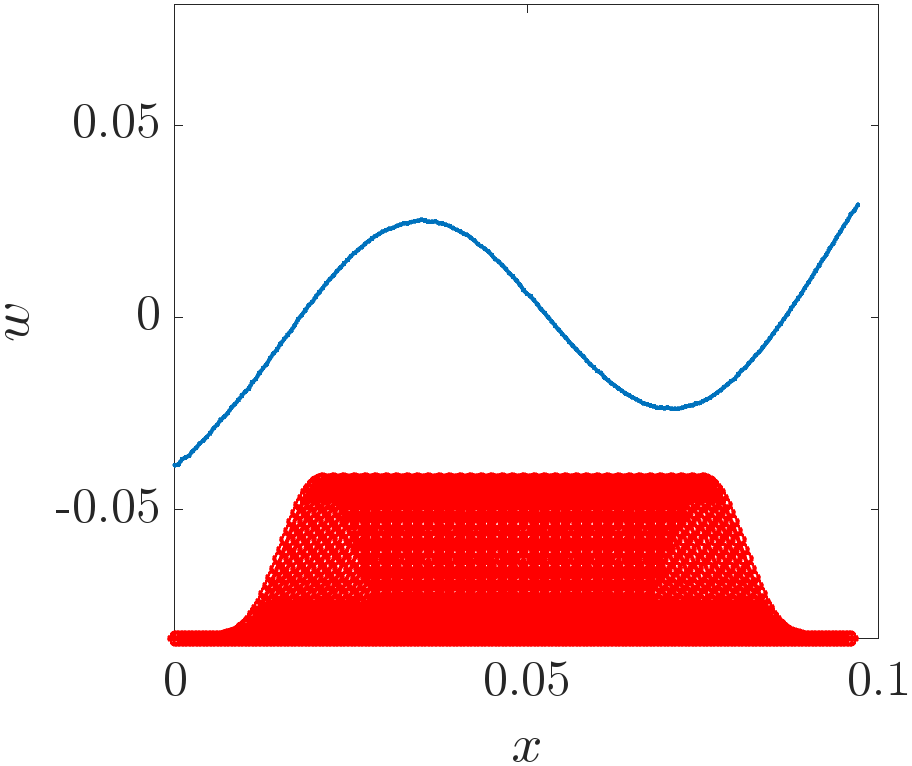}}
        \hfill
        \subfloat[]{\includegraphics[width=0.47\textwidth]{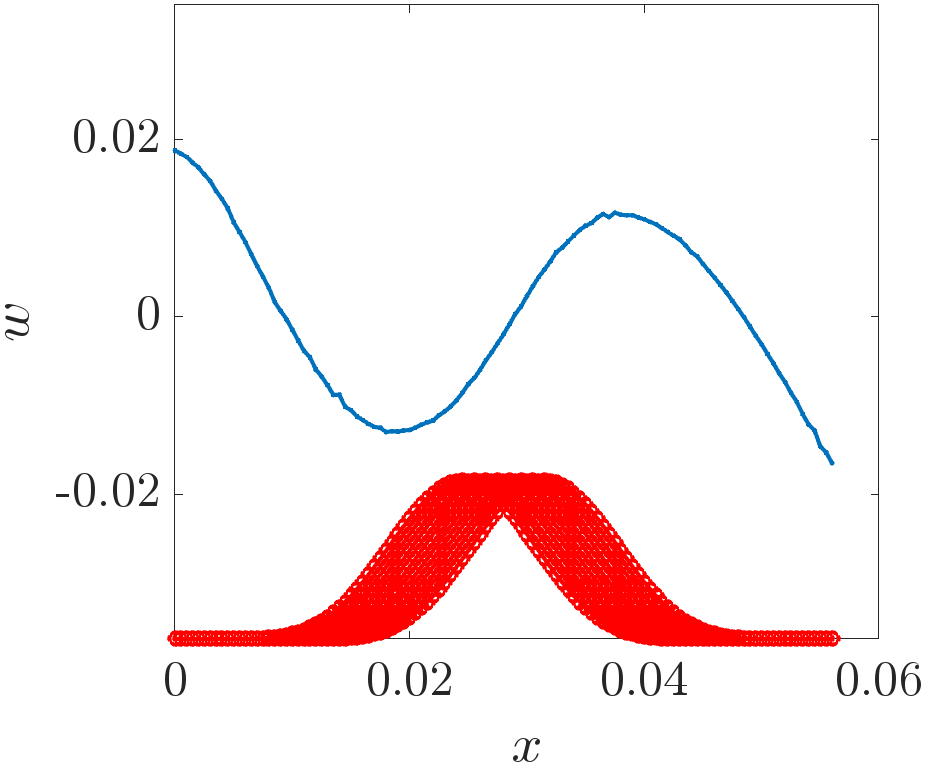}}
        \caption[WSINDy test functions]{Test functions $\phi_1(x)$. Each plot shows the first data point in time (blue) with the spatial portion of the test function (red) over the spatial domain -- Al in (a), IE in (b). The test functions are piecewise polynomials that decay to zero near the boundaries of the domain. The width and order of the test functions are chosen automatically within the WSINDy algorithm (see parameters $m_x, m_t, p_x, p_t$ in Table~\ref{tab:wsindyparams}) and determine how many test functions are used over the spatial domain.}
        \label{fig:wsindy-testfunctions}
    \end{figure}

\begin{figure}[h!]
        \centering
        \subfloat[]{\includegraphics[width=0.47\textwidth]{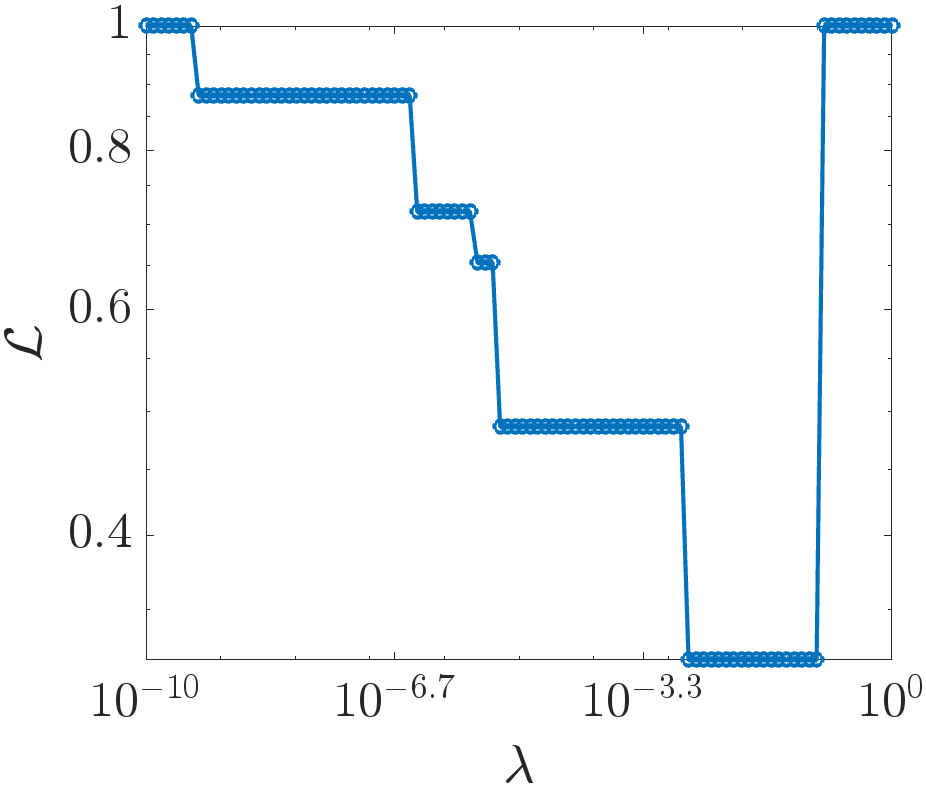}}
        \hfill
        \subfloat[]{\includegraphics[width=0.47\textwidth]{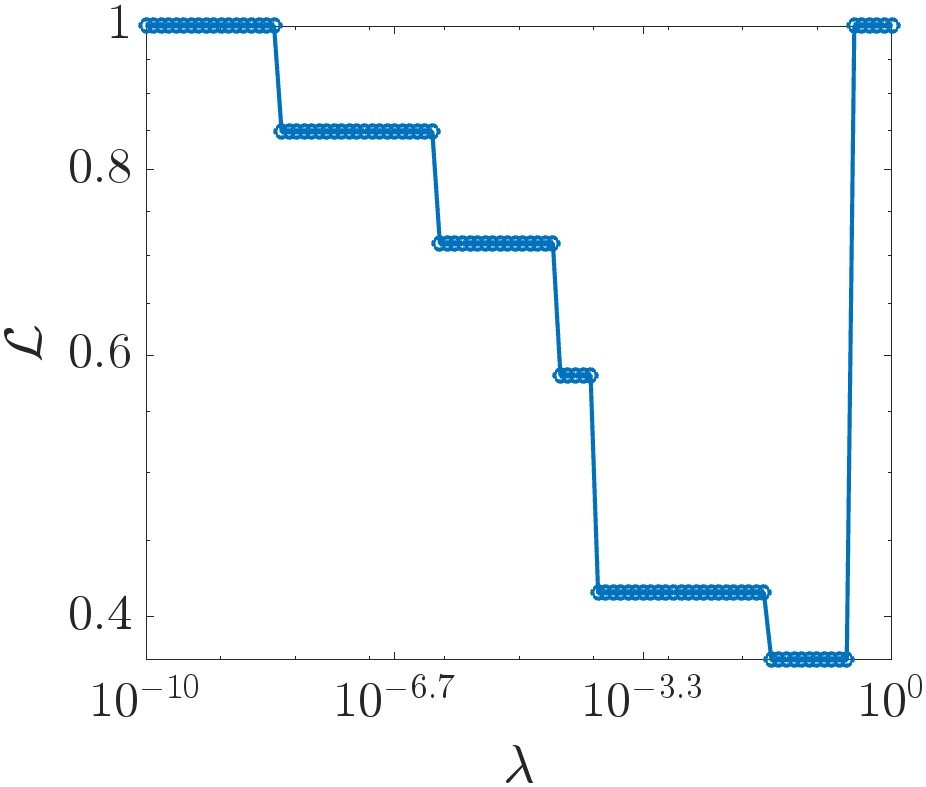}}
        \caption[WSINDy lambda]{Optimization of the sparsity parameter in the modified sequential-thresholding least-squares algorithm. Each plot shows the loss function $
        \mathcal{L}(\lambda) = \frac{\Vert\mathbf{G}(\mathbf{c}_\lambda - \mathbf{c}_{LS})\Vert_2}{\Vert\mathbf{G}\mathbf{c}_{LS}\Vert_2} + \frac{\Vert \mathbf{c}_\lambda \Vert_0}{J}$, 
         vs $\lambda$ for each dataset -- Al in (a) and IE in (b). The optimal value $\hat{\lambda}$ for each data set is the smallest value of $\lambda$ which minimizes the loss function. The optimal $\hat{\lambda}$ values are given in Table~\ref{tab:wsindyparams}.}
        \label{fig:wsindy-lambda}
    \end{figure}

    Table~\ref{tab:wsindyparams} summarizes the hyperparameter values for the analyses in this work. Fig.~\ref{fig:wsindy-testfunctions} shows the spatial component $\phi_1(x)$ of the separable test functions $\psi^{(k)}(x,t)$ for the first time point of each dataset based on the selection of the hyperparameters. Fig.~\ref{fig:wsindy-lambda} shows the optimization of the sparsity parameter $\lambda$. Note these plots are only shown for the non-ensemble equation discovery cases. 

    \begin{table}[h!]
        \centering
        \begin{tabular}{|l|l|l|l|l|}
        \hline
        & \textbf{Al} & \textbf{IE} & \textbf{Ensemble Al} & \textbf{Ensemble IE} \\ \hline
        \hline 
        Data size & $195 \!\times 1501 $ & $133 \!\times 2801$ & $195 \!\times 150 - 1501 $ & $133 \!\times 280 - 2801$ \\ \hline
        $\hat{\tau}$ & $7.3$ & $8.5$ & $7.3 - 7.35$ & $8.3 - 8.55$ \\ \hline
        ($m_x,\ m_t$) & ($42, 82$) & ($49, 84$) & ($42-43, 66 - 82$) & ($49 - 50, 72 - 84$) \\ \hline
        ($s_x,\ s_t$) & ($3, 30$) & ($2, 56$) & ($3, 3 - 30$) & ($2, 5 - 56$) \\ \hline
        ($p_x,\ p_t$) & ($8, 7$) & ($8, 7$) & ($8, 7$) & ($8, 7$) \\ \hline
        $\hat{\lambda}$ & $1.874e-03$ & $2.42e-02$ & $1.177e-03 - 4.7e-03$ & $2.42e-02 - 7.743e-02$ \\ \hline
        ($\gamma_w, \gamma_x, \gamma_t$) & \begin{tabular}[c]{@{}l@{}}($1.00$, \\ $2.42e+02$, \\ $2.85e+05$)\end{tabular} & \begin{tabular}[c]{@{}l@{}}($1.00$, \\ $2.08e+02$, \\ $1.39e+05$)\end{tabular} & \begin{tabular}[c]{@{}l@{}}($1.00$, \\ $2.37e+02 - 2.42e+02$, \\ $3.16e+04 - 2.85e+05$)\end{tabular} & \begin{tabular}[c]{@{}l@{}}($1.00$, \\ $2.04e+02 - 2.12e+02$, \\ $1.54e+04 - 1.39e+05$)\end{tabular} \\ \hline
        Size($\mathbf{G}$) & $1655 \!\times 6 $ & $384 \!\times 6 $ & $37 - 1655 \!\times 6 $ & $182 - 384 \!\times 6 $ \\ \hline
        Condition \#$(\mathbf{G})$ & $1.18e+04$ & $4.44e+04$ & $1.18e+04 - 1.59e+05$ & $4.44e+04 - 1.25e+05$ \\ \hline
        \end{tabular}
        \caption{WSINDy and ensemble WSINDy parameters. In the ensemble columns, fields with ranges of values give the minimum and maximum of the parameter over the ensembles and single values indicate that the parameter was the same across all data ensembles.}
        \label{tab:wsindyparams}
    \end{table}

